\documentclass[11pt]{article}
\usepackage{graphicx}
\usepackage{amsmath,amsfonts,amsthm,amssymb}
\usepackage{latexsym,amsmath}
\usepackage{graphicx,psfrag}

    \psfrag{N1}{$N_1$}
    \psfrag{N2}{$N_2$}
    \psfrag{C1}{$c_1$}
    \psfrag{C2}{$c_2$}
    \psfrag{l1}{$\lambda_1$}
    \psfrag{l2}{$\lambda_2$}
    \psfrag{m1}{$\mu_1$}
    \psfrag{m2}{$\mu_2$}
    \psfrag{n1}{$\mu_1$}
    \psfrag{n2}{$\mu_2$}
    \psfrag{b1}{$\beta_1$}
    \psfrag{b2}{$\beta_2$}
    \psfrag{a}{$a$}
    \psfrag{b}{$b$}
    \psfrag{c}{$c$}
    \psfrag{d}{$d$}
    \psfrag{U}{$\underset{\varphi}\bigcup$}

\newtheorem{theorem}{Theorem}
\newtheorem{lemma}{Lemma}
\newtheorem*{corollary}{Corollary}

\theoremstyle{definition}
\newtheorem*{definition}{Definition}

\theoremstyle{remark}
\newtheorem*{remark}{Remark}
\newtheorem*{ack}{Acknowledgement}

\numberwithin{equation}{section}

\newcommand{\Int}{\operatorname{Int}}
\renewcommand{\Bbb}{\mathbb}

\begin{document}

\title
{3-manifolds that admit knotted solenoids as attractors}
\maketitle
\begin{center}
BOJU JIANG, jiangbj@math.pku.edu.cn
\\
{\it Depart. of Mathematics,
Peking University, Beijing  100871, China\\
}

\end{center}

\begin{center}
YI NI,  yni@princeton.edu\\
{\it Depart. of Mathematics,
Princeton University,  NJ  08544, U. S. A.\\
}
\end{center}

\begin{center}
SHICHENG WANG, wangsc@math.pku.edu.cn\\
{\it Depart. of Mathematics, Peking University, Beijing 100871,
China\\
}
\end{center}

\begin{abstract} Motivated by the study in Morse theory and
Smale's work in dynamics, the following questions are studied and
answered: (1) When does a 3-manifold admit an automorphism having
a knotted Smale solenoid as an attractor? (2) When does a
3-manifold admit an automorphism whose non-wandering set consists
of Smale solenoids? The result  presents some intrinsic symmetries
for a class of 3-manifolds.
\end{abstract}

\maketitle

\section{Introduction}
\label{s:1}

The solenoids are first defined in
mathematics by Vietoris in 1927 for 2-adic case and by others
later in general case, which can be presented either in an
abstract way (inverse limit of self-coverings of circles) or in a
geometric way (nested intersections of solid tori). The solenoids
are introduced into dynamics by Smale as hyperbolic attractors in
his celebrated paper \cite{S}.

Standard notions in dynamics and in 3-manifold topology will be
given in Section \ref{s:2}. The new definitions are the following:

Let $N=S^1\times D^2$, where $S^1$ is the unit circle and $D^2$ is
the unit disc. Both $S^1$ and $D^2$ admit ``linear structures''. Let
$e:N\to N$ be a ``linear'', $D^2$-level-preserving embedding such
that (a) $e(S^1\times *)$ is a $w$-string braid in $N$ for each
$*\in D^2$, where $w>1$ in an integer; (b) for each $\theta\in S^1$,
the radius of $e(\theta \times D^2)$ is $1/w^2$.

\begin{definition}
Let $M$ be a 3-manifold and $f: M\to M$ be a diffeomorphism. If
there is a solid torus $N\subset M$ such that $f|N$ (resp.\
$f^{-1}|N$) conjugates $e:N\to N$ above, we call
$S=\cap_{h=1}^{\infty} f^{h}(N)$ (resp.\ $S=\cap_{h=1}^{\infty}
f^{-h}(N)$) a {\em Smale solenoid}, which is a hyperbolic
attractor (resp.\ repeller, or negative attractor) of $f$, and we
also say $M$ admits $S$ as a Smale solenoid attractor and $N$ is
a defining solid torus of $S$.
\end{definition}

Smale solenoid in the above definition carries more  information
than a solenoid as a topological space. It also carries the
information of braiding of $e(N)$ in $N$ and the knotting and
framing of $N$ in $M$, in addition to  the information that it is
a hyperbolic attractor of a diffeomorphism $f:M\to M$.

Say a Smale solenoid $S\subset M$ is {\em trivial}, if the core of
a defining solid torus $N$ bounds a disc in $M$, otherwise we say
$S$ is {\em knotted}.

\begin{theorem}
\label{t:1}
Suppose $M$ is a closed orientable
3-manifold. There is a diffeomorphism $f: M\to M$ such that the
non-wandering set  $\Omega(f)$ contains a knotted Smale solenoid
\emph{IF} and \emph{ONLY IF} the manifold $M$ has a lens space $L(p,q)$,
with $p\ne 0,\pm1$, as a prime factor.
\end{theorem}

\begin{theorem}
\label{t:2}
Suppose $M$ is a closed orientable
3-manifold. There  is a diffeomorphism $f: M\to M$ with the
non-wandering set $\Omega(f)$ a union of finitely many Smale
solenoids, \emph{IF} and \emph{ONLY IF} the manifold $M$ is a lens space $L(p,q)$,
$p\ne 0$.

Moreover for the IF part,  the $\Omega(f)$ can be chosen to be two
explicit $(p+1)$-adic solenoids, $p +1\ne 0, \pm 1$.
\end{theorem}

\begin{corollary}
The diffeomorphism $f$ constructed in the
\emph{IF} part of Theorem \ref{t:2} is $\Omega$-stable, but is not structurally
stable.
\end{corollary}

\subsection*{Motivations of the results}

\subsubsection* {From Morse theory}  Let $f:M\to R$ be a non-degenerate
Morse function. Then the gradient vector field ${\text {grad} f}$
is a dynamical system on $M$ with hyperbolic $\Omega(\text{grad} f)$.
An important aspect of Morse theory is to
use the global information of the singularities of $f$, or
equivalently, the information of $\Omega(\text{grad} f)$, to
provide topological information of the manifold $M$. The classical
examples are: if $\Omega(\text{grad} f)$ consists of  two points,
then $M$ is the sphere by Reeb in 1952 \cite{R}, and if
$\Omega(\text{grad} f)$ consists of three points, then $M$ is a
projective plane like manifold of dimension 2, 4, 8 or 16 proved
by Eells and Kuiper in 1961 \cite{EK}. The ONLY IF part of Theorems \ref{t:1}
and \ref{t:2} are results of this style.

\subsubsection* {From dynamics of Smale's school} In \cite{S}, for a
diffeomorphism $f:M\to M$, Smale introduced the Axiom A, the
strong transversality condition and the no cycle condition for
$\Omega(f)$. Important results in the dynamics school of Smale are
the equivalences between those conditions and various stabilities.
For an Axiom A system $f$, Smale  proved  (Spectral Decomposition
Theorem)  $\Omega(f)$ can be decomposed into the so-called basic
sets. He posed several types of basic sets: (a) Zero dimensional
ones such  as  isolated points and Smale Horse Shoe; (b) Anosov
maps and maps derived from Anosov; (c) expansive ones such as
Smale solenoids.

All those results and notions need examples to testify. Most known
examples are local. It is natural to ask a global question where
topology and dynamics interact:
 For which manifold $M$, is there an $f:M\to M$ such that all the
basic sets of $\Omega(f)$ belong to a single type above?

There is no restriction when $\Omega(f)$ is zero dimensional. The
answer  to the question for Anosov map was given by Porteous in
1974 \cite{Po}. The ONLY IF part of Theorem \ref{t:2} gives an answer about
Smale solenoids for 3-manifolds. The Corollary also provides
3-dimensional global examples to testify the notions of stability.

We would also like to point out that there are many nice results
on the interplay of topology and dynamics, mostly for flows.
See \cite{F}, \cite{Su} and \cite{T} for examples.

\subsubsection* {Searching symmetries of manifolds with stability.} A
manifold $M$ admitting a dynamics $f$ such that $\Omega(f)$
consists of two hyperbolic attractors presents a symmetry of the
manifold with certain stability. The sphere, the simplest closed
manifold, admits a hyperbolic dynamics $f$ such that $\Omega(f)$
consists of exactly two points, one is a source, and the other is
a sink. The attractors in this example are the simplest in three
senses: (1) The topology of the attractors are trivial, (2) the
embedding of attractors into the manifolds are trivial, (3) the
restriction of the dynamics $f$ on the attractors are trivial. The
IF part of Theorem \ref{t:2} and the Corollary show more manifolds with such
symmetry when we consider more complicated attractors suitably
embedded into the manifolds.

Indeed we believe that many more 3-manifolds admit such symmetries
if we replace the Smale solenoid by its generalization, the
so-called {\em Smale-Williams solenoid} \cite{W} (the name is suggested
in \cite{Pe}).

\subsection*{The structure of the paper}
For the convenience of the
readers from both dynamics and 3-manifold topology, we list the
needed notions and facts in dynamics and in 3-manifold topology in
Section \ref{s:2}.
Sections \ref{s:3}, \ref{s:4} and \ref{s:5} are devoted respectively  to  the
proofs of the ONLY IF parts of  Theorems \ref{t:1} and \ref{t:2}, the IF parts of
Theorems \ref{t:1} and \ref{t:2}, and the Corollary. Most notions  in dynamics
mentioned in Section \ref{s:2} are only used in Section \ref{s:5}. To the authors,
the most interesting part of the paper is the discovery of the IF
part of Theorem \ref{t:2} and its explicit constructive proof. Since such
an explicit constructive proof is difficult to generalize to the
case of Smale-Williams solenoids, we wonder if there is an
alternative proof for the IF part of Theorem \ref{t:2}.

\section{Notions and facts in dynamics and in 3-manifold topology}
\label{s:2}

\subsection*{From Dynamics}

Everything in this part can be found in \cite{Ni}, unless otherwise
indicated.

Assume $f:M\to M$ is a diffeomorphism of a compact $n$-manifold
$M$.

An {\em invariant set} of $f$ is a subset $\Lambda\subset M$ such
that $f(\Lambda)=\Lambda$. A  point $x\in M$ is {\em
non-wandering} if for any neighborhood $U$ of $x$, $f^n(U)\cap
U\ne \emptyset $ for infinitely many integers $n$. Then
$\Omega(f)$, the {\em non-wandering set} of $f$, defined as the
set of all non-wandering points,  is an $f$-invariant closed set.
A set $\Lambda\subset M$ is an {\em attractor} if there exists a
closed neighborhood $U$ of $\Lambda$ such that $f(U)\subset\Int
U$, $\Lambda=\bigcap_{h=1}^\infty f^h(U)$, and
$\Lambda=\Omega(f|U)$.

Say $f$ is {\em structurally stable} if all diffeomorphisms
$C^1$-close to $f$ are conjugate to $f$. Say  $f$ is {\em
$\Omega$-stable} if all diffeomorphisms  $C^1$-close to $f$
preserve the structure of $\Omega(f)$.

A closed invariant set $\Lambda$ of  $f$ is {\em
hyperbolic} if there is a continuous $f$-invariant splitting of
the tangent bundle $TM_\Lambda$ into {\em stable} and {\em
unstable bundles} $E^s_\Lambda\oplus E^u_\Lambda$ with
\begin{alignat*}{3}
\|Df^m(v)\| &\leq C\lambda^{-m}\|v\| \quad && \forall v\in
E^s_\Lambda,\ && \forall m>0,
\\
\|Df^{-m}(v)\| &\leq C\lambda^{-m}\|v\| \quad && \forall v\in
E^u_\Lambda,\ && \forall m>0,
\end{alignat*}
for some fixed $C>0$ and $\lambda>1$.

{\sc The Axiom A.} The diffeomorphism $f:M\to M$ satisfies {\em Axiom A} if
(a) the non-wandering set $\Omega(f)$ is hyperbolic; and
(b) the periodic points of $f$ are dense in $\Omega(f)$.

{\sc Spectral Decomposition Theorem.} For $f:M\to M$ satisfying
Axiom A, $\Omega(f)$ can be decomposed in a unique way into
finitely many disjoint sets $B_1,...,B_k$, so that each $B_i$ is
closed, $f$-invariant and contains a dense $f$-orbit.

The $B_i$ in the decomposition above are usually referred to as
{\em basic sets}.

 {\sc Stable Manifold Theorem.}  Suppose $\Omega(f)$ is
hyperbolic. Then for each $x\in \Omega(f)$, the sets
$W^s(x,f)=\{y\in M|\lim_{j\to \infty}d(f^{j}(y),f^{j}(x))=0\}$ and
$W^u(x,f)=\{y\in M|\lim_{j\to \infty}d(f^{-j}(y),f^{-j}(x))=0\}$
are smooth, injective immersions of the $E_x^s$ and $E_x^u$
respectively. Moreover, they are tangent to $E_x^s$ and $E_x^u$ at
$x$ respectively.

$W^s(x,f)$ and $W^u(x,f)$ in the theorem are known as the {\em
stable and unstable manifold of $f$ at $x$}.

{\sc The Strong Transversality Condition.}  For all $x,y\in \Omega
(f)$, the stable and unstable manifolds, $W^s(x,f)$ and
$W^u(y,f)$, are transverse.

{\sc The No Cycle Condition.} An $n$-cycle of Axiom A system is a
sequence of basic sets $\Omega_0$, $\Omega_1$, ..., $\Omega_n$
with $\Omega_0=\Omega_n$ and $\Omega_i\ne \Omega_j$ otherwise, and
such that $W^u(\Omega_{i-1})\cap W^s(\Omega_i)\ne \emptyset$. An
Axiom A system satisfies the {\em no-cycle condition} if it has no
$n$-cycle for all $n\ge 1$.

{\sc Stability Theorem.} (See the survey paper \cite{Ha})  (a) the
Axiom A and the strong transversality condition of $\Omega(f)$ are
equivalent to the structural stability of $f$.
(b) the Axiom A and the no cycle condition of $\Omega(f)$ are
equivalent to the $\Omega$-stability of $f$.

\subsection*{From 3-manifold theory}

Everything in this part can be found in \cite{He}, unless otherwise
indicated.

Let $M$ be a $3$-manifold and $S$ an embedded $2$-sphere
separating $M$. Let $M_1$ and $M_2$ be the two $3$-manifolds
obtained by splitting $M$ along $S$ and capping-off the two
resulting $2$-sphere boundary components by two $3$-cells. Then
$M$ is {\em a connected sum} of $M_1$ and $M_2$, written
$M_1\#M_2$.

A $3$-manifold $M\ne S^3$ is {\em prime} if $M=M_1\#M_2$ implies
one of $M_1$, $M_2$ is $S^3$.

Let $F$ be a connected compact 2-sided surface properly embedded
in $M$. $F$ is said to be {\em compressible} if either $F$ bounds
a 3-ball, or there is an essential, simple closed curve on $F$
which bounds a disk in $M$; otherwise, $F$ is said to be {\em
incompressible}.

The following three  results in 3-manifold topology are
fundamental.

{\sc Kneser-Milnor's Prime Decomposition Theorem.} Every closed
orientable $3$-manifold $M\ne S^3$ can be expressed as a connected
sum of a finite number of prime factors. Furthermore, the
decomposition is unique up to order and homeomorphism.

{\sc Haken's Finiteness Theorem.} Let $M$ be a compact orientable
3-manifold. Then  the maximum number of pairwise disjoint,
non-parallel, closed connected incompressible surfaces in $M$,
denoted by $h(M)$, is a finite integer $\ge 0$.

{\sc Papakyriakopoulos's Loop  Theorem.} Let $M$ be a compact
orientable 3-manifold and $S\subset M$ a closed  orientable
surface. If the homomorphism $i_*: \pi_1(S)\to \pi_1(M)$ induced
by the embedding $i: S\to M$ is not injective, then there is an
embedded disc $D\subset M$ such that $D\cap S=\partial D$ and
$\partial D$ is an essential circle in $S$.

For the definition of the {\em lens space $L(p,q)$}, see Section
\ref{s:4}.

\section{Proof of the ONLY IF parts of Theorems \ref{t:1} and \ref{t:2}}
\label{s:3}

We first prove the ONLY IF part of Theorem \ref{t:1}.

\begin{proof}
Suppose $f: M\to M$ has a knotted Smale solenoid $S$
as an attractor. Then $S=\cap_{h=1}^{\infty} f^{h}(N)$, and
$\overline{M-N}\subset {M-f(N)}$, where $N$  is a defining solid
torus of $S$.

Since $f$ is a global homeomorphism, $\overline{M-N}$ and
$\overline{M-f(N)}$ are homeomorphic.

Suppose first that $\partial\overline {M-N}$ is an incompressible
surface in $\overline{M-N}$. By Haken's Finiteness Theorem,
$h(M)$, the maximum number of pairwise disjoint, non-parallel,
closed incompressible surfaces in $M$, is a finite integer. Since
the winding number $w$ of $f(N)$ in $N$ is $>1$,
$\partial\overline{M-N}$ is incompressible in $\overline{N-f(N)}$
and is not  parallel to $\partial \overline{M-f(N)}$. It follows
that for any set $F$ of disjoint, non-parallel, incompressible
surfaces of $\overline{M-N}$, $\partial \overline{M-f(N)}\cup F$
is a set of disjoint non-parallel closed incompressible surfaces
in $\overline{M-f(N)}$. Hence $h(\overline{M-f(N)})$ is larger
than $h(\overline{M-N})$, which contradicts the fact that
$\overline{M-N}$ and $\overline{M-f(N)}$ are homeomorphic.

By the last paragraph, $\partial\overline {M-N}$ is compressible
in $\overline{M-N}$. This means  there is a properly embedded disc
$(D, \partial D)\subset (\overline{M-N},\partial \overline{M-N})$
such that $\partial D$ is an essential circle in $\partial N$.
Cutting $\overline{M-N}$ along $D$, we get a 3-manifold, denoted
by $M_1$, with $\partial M_1$ a 2-sphere containing two copies
$D_1$ and $D_2$ of $D$. Let $S_*$ be a boundary parallel 2-sphere
in the interior of $M_1$. Now identifying $D_1$ and $D_2$, we get
back to $\overline{M-N}$  and $S_*$ separates a punctured solid
torus from $\overline{M-N}$; finally we glue back $N$ with
$\overline {M-N}$ to get $M$, $S_*$ separates a punctured lens
space from $M$, i.e., $M$ contains a lens space $L$ as a prime
factor.

If $L$ is  $S^3$, then it  is easy to see the core of $N$ bounds a
disc, which contradicts the assumption that $S$ is knotted. If
$L=S^2\times S^1$, then $N$ carries a generator $\alpha$ of
$\pi_1(S^2\times S^1)=Z$. Since $f(N)$ is a $w$-string braid in
$N$, we have $f_*(\alpha)=w\alpha$. Since $f$ is a
homeomorphism, $f_*$ is an isomorphism. Hence $w=1$, and we reach
a contradiction.

We have finished the proof of the ONLY IF part of Theorem \ref{t:1}.
\end{proof}

We are going to prove the ONLY IF part of Theorem \ref{t:2}.

Suppose $\Omega(f)$ is a union of Smale solenoids $S_1, ... ,
S_n$. Then for each $i=1,...,n$, it is known (more or less
directly from the definition) that

(i) $f|S_i$ is hyperbolic and the periodic points of $f$ are dense
in $S_i$;

(ii) $S_i$ is an $f$-invariant closed set and there is a dense
$f$-orbit in $S_i$.

Then $f$ satisfies the Axiom A by (i). By Spectral Decomposition
Theorem, $\Omega(f)$ can be decomposed in a unique way into
finitely many disjoint basic sets $B_1,...,B_k$, so that each
$B_i$ is closed, $f$-invariant and contains a dense $f$-orbit.

By (ii), each $S_i\subset B_l$ for some $l=1,...,k$. Then from the
facts that $S_i$ is an attractor of $f$ (or of $f^{-1}$) and that
$B_l$ contains a dense $f$-orbit,  there is a point $x\in \Int
U_i$ so that its $f$-orbit $o(x)$ is dense in $B_l$, where $U_i$
is a closed neighborhood of $S_i$ mentioned in the definition of
an attractor. Then it is clear that $x\in \Omega(f|U_i)$. Hence
$x\in S_i$, thus $B_l=\overline{o(x)}\subset S_i$, so we must have
$S_i=B_l$. Hence each $S_i$ is a basic set of $\Omega(f)$ and in
particular, $\Omega(f)$ is a disjoint union of finitely many Smale
solenoids.

Now the ONLY IF part of Theorem \ref{t:2} follows from the following Lemma
\ref{l:1} and Lemma \ref{l:2}.

\begin{lemma}
\label{l:1}
Suppose $f:M\to M$ is a diffeomorphism and
$\Omega(f)$ is a disjoint union of finitely many Smale solenoids.
Then $\Omega(f)$ is a union of two solenoids, one is an attractor
of $f$ and the other is an attractor of $f^{-1}$.
\end{lemma}

\begin{proof}
Suppose
\begin{equation}
\Omega(f)=S_1\cup S_2\cup ...\cup S_n \qquad n \ge 1, \tag{1}
\end{equation}
is a disjoint union of solenoids, where either
$S_i=\cap_{h=1}^{\infty} f^{h}(N_i)$ if $S_i$ is an attractor of
$f$ or $S_i=\cap_{h=1}^{\infty} f^{-h}(N_i)$ if $S_i$ is an
attractor of $f^{-1}$. Without loss of generality, we assume that
the $N_i$'s have been chosen so that  $N_i\cap N_j=\emptyset$ if
$i\ne j$ (since $S_i\cap S_j=\emptyset$ for $i\ne j$), and some
$S_i$ is an attractor of $f$ (otherwise replace $f$ by $f^{-1}$).
Then by re-indexing if necessary we assume that $S_1, ..., S_k$
are attractors of $f$ and the remaining $S_j$ are attractors of
$f^{-1}$; henceforth, $k$ is the number of attracting solenoids.
So we can assume
\begin{equation}
f(N_i) \subset \Int N_i, \qquad i=1,...,k\le n,\tag {2}
\end{equation}
and
\begin{equation}
f^{-1}(N_j)\subset \Int N_j, \qquad j=k+1,...,n.\tag {3}
\end{equation}

For $i=1,\dots,k$, let $V_i= \cup_{h=1}^{\infty} f^{-h}(\Int N_i)$.
Since $f$ is a homeomorphism, $V_i$ is open. Moreover
\begin{equation}
\qquad V_i\cap V_j =\emptyset, \qquad 1\le i<j\le k\tag {4}
\end{equation}
and
\begin{equation}
f(V_i)=V_i, \tag {5}
\end{equation}
by the assumptions (2) and $N_i\cap
N_j=\emptyset$ for $i\ne j$.

First we suppose $n=1$. Now Let $Y_1= M - \Int N_1$. Then $Y_1$ is
compact and $f^{-1}(Y_1)\subset Y_1$ by (2). Therefore
$\Omega(f)=\Omega (f^{-1})$ intersects $Y_1$, which contradicts
(1).

So $n>1$. Suppose $k>1$. Let $Y_2=M-\cup_ {j=k+1}^n S_j$.
For each $i=1,\dots,k$, $V_i\subset Y_2$.
$Y_2$ is connected, so it cannot be a disjoint
union of $k>1$ open sets. Hence $Y_3=M-((\cup_{j=k+1}^n S_j)\cup
(\cup_{i=1}^k V_i))$ is not empty. Suppose $x\in Y_3$, since $S_j$
is compact, we can choose $N_j$ sufficiently small in order that
$x\notin \cup_{j=k+1}^n \Int N_j$.  Then
$Y_4=M-((\cup_{j=k+1}^n \Int N_j)\cup (\cup_{i=1}^k V_i))$ is
compact and is not empty. By (3),(5), we have $f(Y_4)\subset Y_4$.
Hence $\Omega (f)\cap Y_4 \ne \emptyset$, a contradiction.

We have proved that $f$ has exactly one attractor. By the same
reason $f^{-1}$ also has exactly one attractor, therefore $n=2$
and the lemma is proved.
\end{proof}

\begin{lemma}
\label{l:2}
Let $M$ be a closed orientable 3-manifolds. If
$f: M\to M$ is a diffeomorphism with $\Omega(f)$ a union of two
disjoint Smale solenoids, then $M$ is a lens space and $M$ is not
$S^2\times S^1$.
\end{lemma}

\begin{proof}
Suppose $\Omega(f)$ is a union of two disjoint
solenoids $S_1$ and $S_2$. We may further assume that
\begin{equation}
S_1=\cap_{h=1}^{\infty} f^{h}(N_1), \qquad S_2=\cap_{h=1}^{\infty}
f^{-h}(N_2), \qquad N_1\cap N_2=\emptyset.\tag {6}
\end{equation}
We have $\cup_{h=1}^{\infty} f^{-h}(\Int N_1)=M-S_2$.
It follows that
\begin{equation}
f^n(\partial N_2)\subset \Int N_1, \quad M-N_1\subset f^n(N_2)
\quad \text{ for some large integer $n>1$}.\tag {7}
\end{equation}

Since $H_2(N_1,Z)=0$, $\partial f^n(N_2)$ separates $N_1$ into two
parts $Y'$ and $Y''$ with $\partial Y'=\partial f^n(N_2)$ and
$\partial Y''$ has two components.

The homomorphism $i_*: \pi_1(\partial f^n(N_2))\to \pi_1(N_1)$
induced by the embedding $i: \partial f^n(N_2)\to N_1$ is not
injective, since $\pi_1(N_1)=Z$ and $\pi_1(f^n(\partial
N_2))=Z\oplus Z$. By  the Loop Theorem, $\partial f^n( N_2)$ is
compressible in $N_1$, that is, there is an embedded disc
$D\subset N_1$ such that $D\cap
\partial N_2=\partial D$ and $\partial D$ is an essential circle
in $\partial f^n(N_2)$. Since the solid torus $N_1$ is
irreducible, a standard argument shows that $\partial f^n(N_2)$
bounds a solid torus $N'$ in $N_1$, and therefore we have $N'=Y'$.
Then by (7), we have
\begin{equation}
M=(M-N_1)\cup_{\partial N_1} N_1= M-N_1\cup_{\partial N_1}
Y''\cup_{\partial f^n(N_2)} Y'= f^n(N_2)\cup_{\partial f^n(N_2)}
N'. \tag {8}
\end{equation}
Hence $M$ is obtained by identifying  two solid tori $f^n(N_2)$
and $N'$ along their common boundary. So $M$ is a lens space.

Since $f$ is a homeomorphism, $N''=f^{-n}(N')$ is also a solid
torus and $M$ is obtained by identifying  two solid tori $N_2$
and $N''$ along their boundary. Now $f^{-1}(N_2)$ is a $w$-string
braid in $N_2$, $w>1$. That  $M$ is not $S^2\times S^1$ can be
proved as before.
\end{proof}

\section{Proof of the IF parts of Theorems \ref{t:1} and \ref{t:2}}
\label{s:4}


\begin{figure}
    \includegraphics[width=\textwidth]{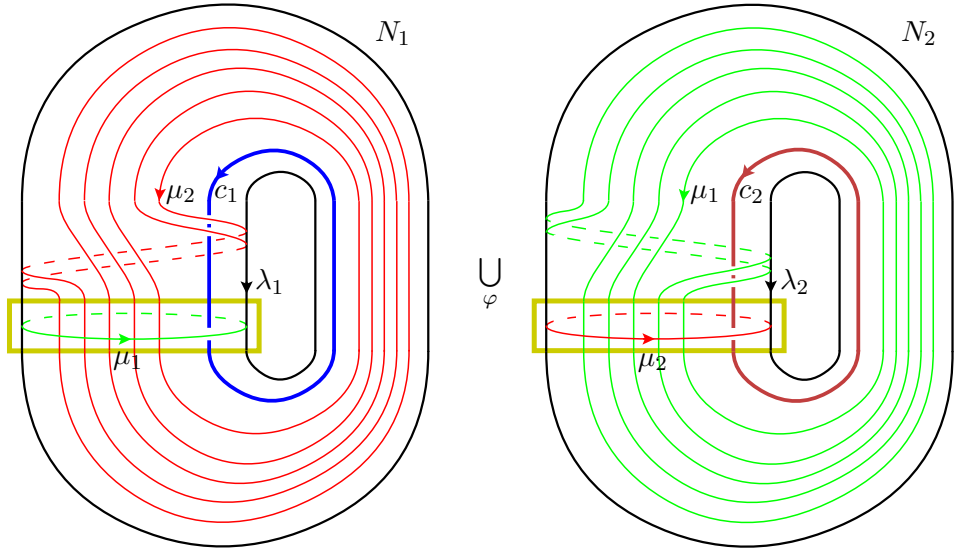}
    \caption{Lens space $L(p,q)$ as union of solid tori $N_1\underset{\varphi}\cup N_2$}
\end{figure}

Suppose $M$ is the lens space $L(p,q)$, where $p>0$ and
$\gcd(p,q)=1$. Then $M$ is the union of two solid tori,
$M=N_1\cup_\varphi N_2$, where the gluing map $\varphi: \partial
N_2 \to\partial N_1$ is an orientation reversing homeomorphism. On
each torus $\partial N_i$, pick a meridian-longitude pair, denoted
$\{\mu_i,\lambda_i\}$, as basis of $H_1(\partial N_i)$. In
$\partial N_1$, $\varphi(\mu_2)$ is the $(p,q)$-curve, that is
$\varphi(\mu_2)=p\lambda_1+q\mu_1$, while
$\varphi(\lambda_2)=r\lambda_1+s\mu_1$, with $ps-qr=1$. It is
clear that in $\partial N_1$ we have
$\varphi^{-1}(\mu_1)=p\lambda_2-r\mu_2$ and
$\varphi^{-1}(\lambda_1)=-q\lambda_2+s\mu_2$. In Figure~1 , the
case of $M=L(5,2)$ and $\begin{pmatrix} p & q \\ r & s
\end{pmatrix} =\begin{pmatrix} 5 & -2 \\ -2 & 1 \end{pmatrix}$ is
shown as a concrete example.

\subsection*{Proof of the  IF part of Theorem \ref{t:2}}
The IF part of Theorem \ref{t:2}
is equivalent to the following

{\bf Claim.} {\em Suppose $M$ is a lens space $L(p,q)$, $p> 0$.
Then there is a diffeomorphism $f: M\to M$ with $\Omega(f)$ a
union of two $(p+1)$-adic solenoids, one is an attractor, the
other is a repeller.}

We are going to prove this Claim.


\begin{figure}
    \includegraphics[width=\textwidth]{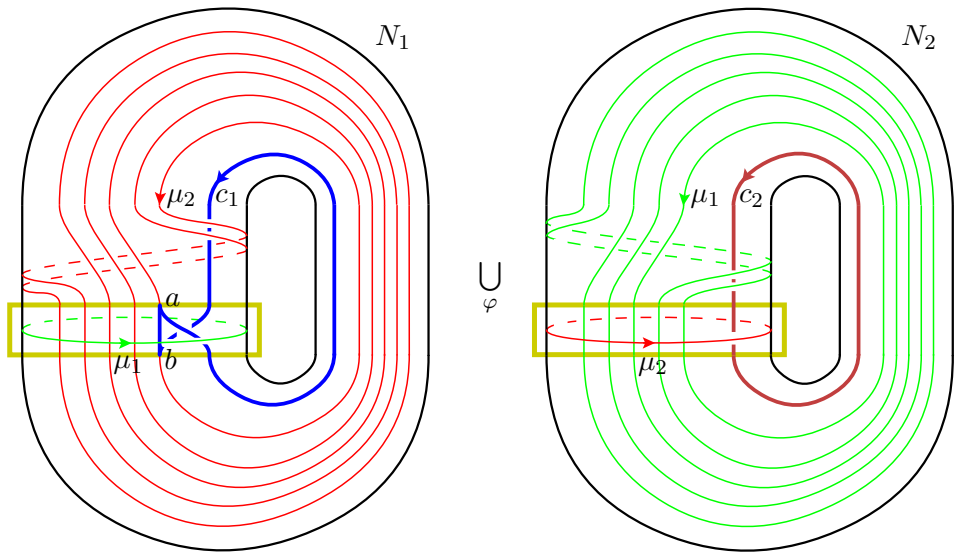}
    \caption{Writhe the core $c_1$ in $N_1$}
\end{figure}

\begin{figure}
    \includegraphics[width=\textwidth]{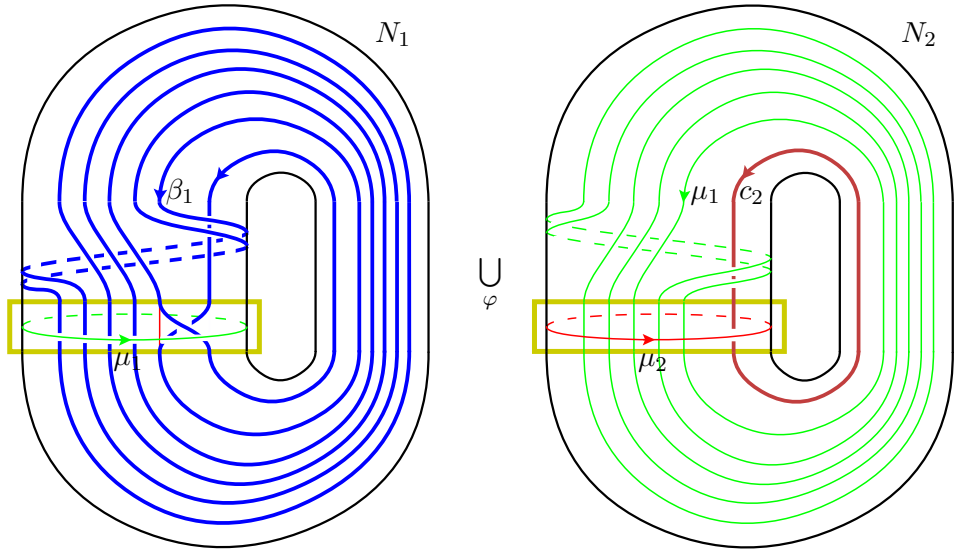}
    \caption{Closed braid $\beta_1$ in $N_1$ and core $c_2$ of $N_2$}
\end{figure}

\begin{figure}
    \includegraphics[width=\textwidth]{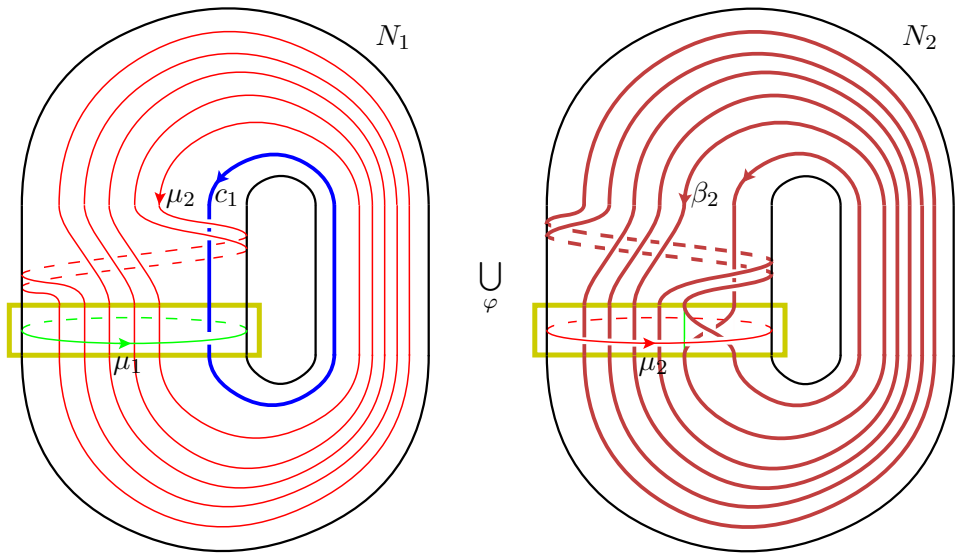}
    \caption{Core $c_1$ of $N_1$ and closed braid $\beta_2$ in $N_2$}
\end{figure}

Denote the oriented cores of $N_1,N_2$ by $c_1,c_2$ ($c_i$ is
homologous to $\lambda_i$ in $N_i$) respectively.
  We do the following
operations to $c_1$, as indicated in Figure~2: Writhe $c_1$
locally, moving a subarc $\overline{ab}$ toward $\partial N_1$ and
identify it with a subarc of $\varphi(\mu_2)$. Since $\mu_2$
bounds a meridian disk in $N_2$, we can push $\overline{ab}$
across the disk, the effect seen in $N_1$ is to replace
$\overline{ab}$ with its complement in $\varphi(\mu_2)$, see
Figure~2. Finally, pushing the obtained curve into $\Int N_1$, we
get a closed braid $\beta_1$ in $N_1$, as indicated in Figure~3.
(In fact, $\beta_1$  is the ``connected sum'' of the ``writhed''
$c^{-1}$ with $\varphi(\mu_2)$ in $N_1$.)
  Do similar operations to $c_2$ in $N_2$, we get a closed braid
$\beta_2$ in $N_2$, as indicated in Figure~4. Now $\beta_1\sqcup
c_2$ and $c_1\sqcup\beta_2$ are two links in $M$.

\begin{lemma}
\label{l:3}
The two links $\beta_1\sqcup c_2$ and
$c_1\sqcup\beta_2$ are isotopic in $M$.
\end{lemma}

\begin{proof}
Recall that $\beta_1$ is obtained by isotoping $c_1$,
thus if we perform the inverse of the above isotopy, we can
transform $\beta_1$ into $c_1$. We will show that the same isotopy
also transforms $c_2$ into $\beta_2$.

From now on, we only use local pictures (represented in the
rectangular frame in Figures~1--4) to show changes in both $N_1$
and $N_2$ simultaneously. The initial local picture of
$\beta_1\sqcup c_2$ is shown in Figure~5-1. In $N_1$ (on the
left), $\beta_1$ is a closed braid in $N_1$, and a segment of
$c_2$ is shown outside of $\partial N_1$. On the right, most
of $\beta_1$ coincides with $\mu_2$, along with the part
slightly outside $\partial N_2$, and $c_2$ is the core of $N_2$.
Our isotopy consists of the following three steps:

{\sc Step 1.} $\mu_2$ bounds a meridian disk in $N_2$, so we can pull $\beta_1$
across the disk. At the same time, a subarc of $c_2$ is pulled into
$N_1$, as indicated in Figure~5-2.

{\sc Step 2.} In the local picture Figure~5-2, $\beta_1$ has a self-crossing.
A local half twist will eliminate this self-crossing, as indicated in Figure~5-3.
  Take care so that a subarc $\overline{cd}$ of $c_2$ lies on
$\varphi^{-1}(\mu_1)$.
  Now compare Figure~5-2 and Figure~5-3, we find an
interesting fact: except for the colors and labels, the left/right part
of Figure~5-2 is the same as the right/left part of Figure~5-3.
This symmetry suggests that the next step is a kind of inverse to Step 1.

{\sc Step 3.} Push the subarc $\overline{cd}$ across the meridian
disk of $N_1$, as indicated in Figure~5-4. We see that $\beta_1$
is deformed to $c_1$, and $c_2$ is deformed to $\beta_2$.
\end{proof}


\begin{figure}
    \includegraphics[width=\textwidth]{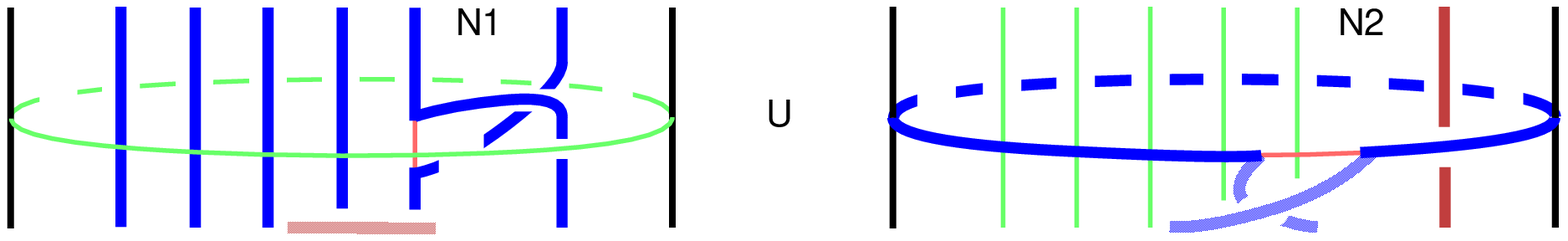}
    \[\Downarrow\]

   \centerline{ Pull across the meridian disk of $N_2$}

    \includegraphics[width=\textwidth]{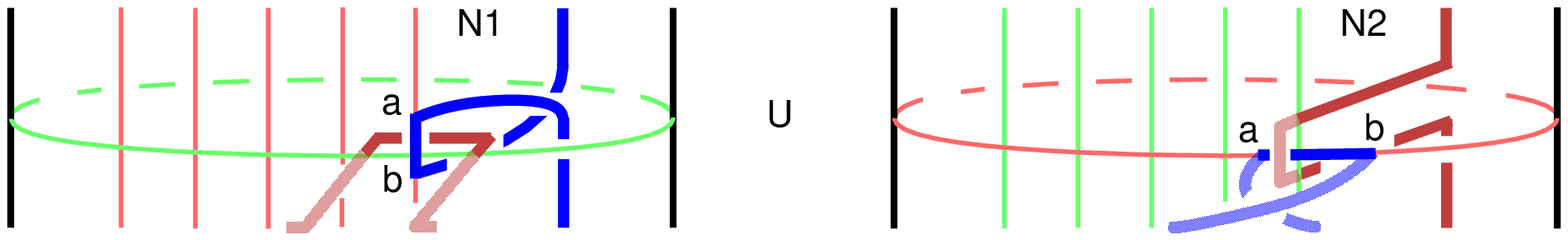}
    \[\Downarrow\]

    \centerline{A local half twist}

    \includegraphics[width=\textwidth]{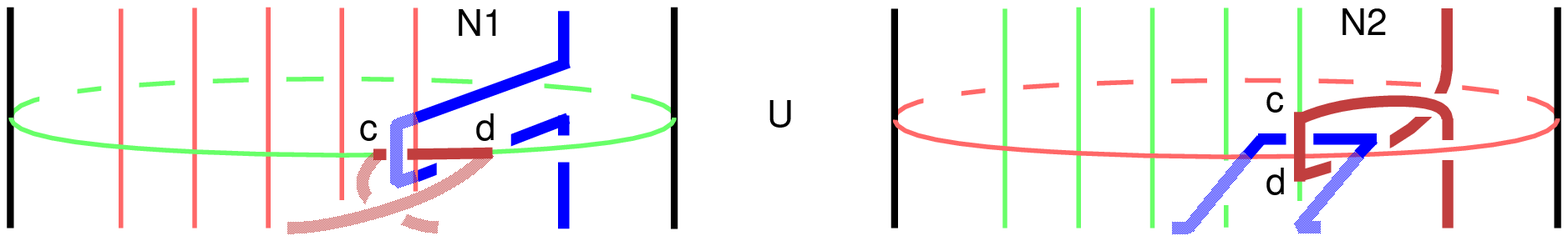}
    \[\Downarrow\]

    \centerline{Push across the meridian disk of $N_1$}

    \includegraphics[width=\textwidth]{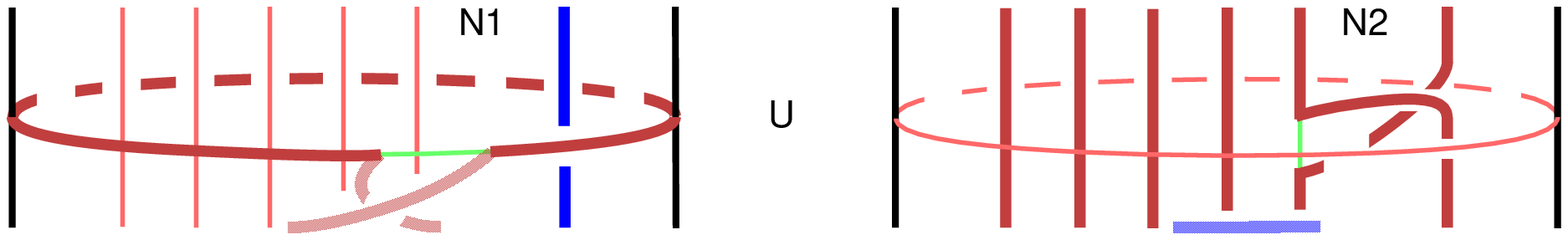}

    \caption{Local pictures of the 3-step isotopy
from $\beta_1\sqcup c_2$ to $c_1\sqcup\beta_2$}
\end{figure}

\begin{proof}[Proof of the  Claim]
Note that $\beta_i$ is a $(p+1)$-string
braid in $N_i=S^1\times D^2_i$. Isotope $\beta_i$ in $N_i$ to meet
all fiber discs $*\times D^2_i$ transversely. Let $\mathcal
N(\gamma)$ denote the closed tubular neighborhood of a closed
curve $\gamma$, and think of $N_i$ as $\mathcal N(c_i)$.

Choose  $\mathcal N(\beta_i)$ to be a disc bundle over $\beta_i$
embedded into $N_i$ so that each disc fiber $\subset$ $* \times
D^2$ (for $*\in S^1$) and has diameter $ < 1/ (p+1)^2$. Moreover
we may assume that $N(\beta_i)$ misses the core $c_i$.

The isotopy provided by Lemma \ref{l:3} that sends $\beta_1\sqcup c_2$ to
$c_1\sqcup\beta_2$ can be adjusted to send $\mathcal
N(\beta_1)\sqcup N_2$ to $N_1\sqcup\mathcal N(\beta_2)$, and to be
``linear'' and ``disc-fiber preserving'' on $\mathcal
N(\beta_1)\sqcup N_2$. Then extend it to a diffeomorphism $f: M\to
M$ which sends $\overline {N_1-N(\beta_1)}$ to $\overline
{N_2-N(\beta_2)}$.

Now the $(p+1)$-adic solenoids $S_1=\cap_{h=1}^{\infty}
f^{-h}(N_1)$ and $S_2=\cap_{h=1}^{\infty} f^{h}(N_2)$ are the
repeller and the attractor of $f$ respectively. Moreover
 for each $x\notin S_1\cup S_2$,
$f^n(x)$ approaches to $S_2$ as $n$ approaches to infinity, hence
$\Omega(f)=S_1\cup S_2$.

We have finished the proof of the Claim, therefore the IF part of
Theorem \ref{t:2}.
\end{proof}

\begin{remark} By repeating the operations in the proof, we see
that in the IF part of Theorem \ref{t:2}, the $\Omega(f)$ can be chosen
to be two explicit $(mp+1)$-adic solenoids, $mp +1\ne 0,\pm 1$.
\end{remark}

\subsection*{Proof of the  IF part of Theorem \ref{t:1}}
Suppose $M=N\#L(p,q)$.
It is  easy to see that the isotopy above that sends
$\beta_1\sqcup c_2$ to $c_1\sqcup\beta_2$ can be adjusted to send
$ N_2$ to $\mathcal N(\beta_2)$, to be ``linear'' and ``disc-fiber
preserving'' on $N_2$, and to be the  identity on a 3-ball $B^3$ in
$N_1$. Therefore there is a diffeomorphism on the $L(p,q) -
\text{int} B^3$ which has a knotted solenoid as a hyperbolic
attractor and is identity on its 2-sphere boundary. Such a
diffeomorphism can be extended to $M$ by the identity on the
punctured $N$.

We have proved the IF part of Theorem \ref{t:1}.

\section{Proof of the Corollary}
\label{s:5}

We start from the end of the proof of the IF part of Theorem \ref{t:2}.

Since $\Omega(f)$ consists of two Smale solenoids $S_1$  and
$S_2$, $\Omega(f)$ meets Axiom A.

To prove the Corollary,   we need the following explicit description
of stable and unstable manifolds of $\Omega(f)$.

First, $S_1$ is the union of stable manifolds of points in $S_1$,
and $S_2$ is the union of  unstable manifolds of points in $S_2$.
Moreover, since $f^{-1}|N_1$ (resp.\ $f|N_2$) preserves the disc
fibers of $N_1$ (resp.\ $N_2$), $\Bbb F_1 =\cup f^{n}(S^1\times
D_1)$ (resp.\ $\Bbb F_2 =\cup f^{-n}(S^1\times D_2)$) provides an
$R^2$-foliation of $L(p,q)-S_2$ (resp.\ $R^2$-foliation of
$L(p,q)-S_1$), which is the union of unstable manifolds of points
in $S_1$ (resp.\ the union of stable manifolds of points in $S_2$).
Hence we have
$$W^s(S_1)=S_1,\quad W^u(S_1)=\Bbb F_1,\quad W^u(S_2)=S_2,\quad W^s(S_2)=\Bbb F_2.$$
Therefore
$$W^s(S_1)\cap W^u(S_2)=\emptyset, \qquad W^s(S_2)\cap W^u(S_1)\ne \emptyset.$$

It is clear that $f$ meets the no cycle condition. Hence $f$ is
$\Omega$-stable by (b) of the Stability Theorem.

This $f$ is not structurally stable by the  Stability Theorem (a)
and the following Lemma \ref{l:4}.

\begin{lemma}
\label{l:4}
$\Bbb F_1$ and $\Bbb F_2$
do not meet transversely.
\end{lemma}

\begin{proof} We need only to prove that $\Bbb F_1|$ and $\Bbb
F_2|$, the restrictions of $\Bbb F_1$ and $\Bbb F_2$ on $\overline
{N_1-N(\beta_1)}$ respectively, do not meet transversely.

Note that $\overline {N_1-N(\beta_1)}$  has two different
$(p+1)$-punctured disc bundle structures  provided by $\Bbb F_1|$
and $\Bbb F_2|$. (An $n$-punctured disc is obtained from the
2-sphere by removing the interior of $n+1$ disjoint sub-discs.)
More directly, one $(p+1)$-punctured disc bundle structure is
induced from the pair ($N_1$, $N(\beta_1)$) and the other is
induced from the pair $(\overline {N_1-N(\beta_1)}\cup N_2,
N_2)=(f^{-1}(N_2), f^{-1}(N(\beta_2))\cong (N_2, N(\beta_2))$.

It is easy to see that the restrictions of two fibrations $\Bbb
F_1|$ and $\Bbb F_2|$ on  $\overline {N_1-N(\beta_1)}$ meet
transversely on $\partial \overline {N_1-N(\beta_1)}$.

Let $F_1$ be a fiber of $\Bbb F_1|$, which is a $(p+1)$-punctured
disc. Suppose $\Bbb F_1|$ and $\Bbb F_2|$  meet transversely on
$\overline {N_1-N(\beta_1)}$. Then the intersections of $F_1$ and
$\Bbb F_2|$ provide a codimension one foliation on $F_1$ which
meets $\partial F_1$ transversely. Now the genus $(p+1)$ closed
surface $D(F_1)$, the double of $F_1$, will admit a codimension
one foliation, which is impossible since $|p+1|>1$.
\end{proof}

We have completed the Proof of the Corollary.

\begin{ack} We would like to thank the referee for his
suggestions which enhanced the paper. We also thank Professors H.\ Duan,
S.\ Gan and L.\ Wen for helpful conversations.
\end{ack}

\bibliographystyle{amsalpha}

\end{document}